\documentclass[11pt, reqno]{amsart}
\usepackage{latexsym}
\usepackage{amsmath, amssymb, amsfonts}
\author{Chengyan Liu}
\address{7003 Brocton Ct., Virginia, VA 22150}
\email{geismn@bigfoot.com}
\title{Riemann Hypothesis}
\begin{document}
\theoremstyle{plain}
\newtheorem{rmmin}{Theorem}
\newtheorem{comin}{Corollary}
\newtheorem{comoo}[comin]{Corollary}
\newtheorem{lfmin}[rmmin]{Theorem}
\newtheorem{hrlem}{Lemma}
\newtheorem{eulem}[hrlem]{Lemma}
\newtheorem{palem}[hrlem]{Lemma}
\newtheorem{tilem}[hrlem]{Lemma}
\newtheorem{volem}[hrlem]{Lemma}
\newtheorem{dedmi}[rmmin]{Theorem}
\newtheorem{hekmi}[rmmin]{Theorem}
\newtheorem{dedme}[hrlem]{Lemma}
\newtheorem{dedfa}[hrlem]{Lemma}
\newtheorem{hkfoo}[hrlem]{Lemma}
\theoremstyle{remark}
\newtheorem*{ounga}{Condition N}
\newtheorem{ovlap}{Observation}
\def\theovlap{\alph{ovlap}}
\theovlap
\newtheorem{totrm}[ovlap]{Observation}
\newtheorem{siuon}{Situation}[ovlap]
\newtheorem{siund}[siuon]{Situation}
\newtheorem{siurd}[siuon]{Situation}
\newtheorem{kisei}{Case}[ovlap]
\newtheorem{kisbe}[kisei]{Case}
\newtheorem{kbead}[kisei]{Case}
\begin{abstract}
Through an equivalent condition on the Farey series set forth by Franel~\cite{frj:fnp} 
and Landau~\cite{lae:bhf}, we prove Riemann Hypothesis for the Riemann 
{zeta-function} $\zeta(s)$ and the Dirichlet $L$-function $L(s, \chi)$. 
\end{abstract}
\maketitle
In the memoir \cite{rib:apg}, Riemann set a stage for investigating the distribution of prime 
numbers through zeta function   
\begin{equation*}\zeta(s) = \sum_{n = 1}^{\infty}{1 \over n^{s}}~,\end{equation*}
where $s$ a complex number, $n$ an integer, $Re(s) > 1$. Augmenting $\zeta(s)$ to the 
whole plane by analytic continuation, Riemann proved functional equation
\begin{equation}\label{fuq:zte}
\pi^{- {1 \over 2} s} \Gamma \biggl({s \over 2} \biggr) \zeta(s) = 
                                  \pi^{- {1 \over 2} (1 - s)} \Gamma \biggl({{1 - s} \over 2} \biggr) \zeta(1 - s)~,
\end{equation}
where $\Gamma(s) = \int_{0}^{\infty} e^{- t} t^{s - 1} dt$ the $\Gamma$ function. Through 
equation \eqref{fuq:zte}, we see that $\zeta(s)$ has trivial zeros at the poles 
of $\Gamma(s/ 2)$, $s = -2, -4, -6, \ldots$. He further suggested the Riemann Hypothesis 
that all nontrivial zeros of $\zeta(s)$ lie on line $Re(s) = 1/ 2$.

One of its equivalent condition is related to the Farey sequences. Denote $r_{\nu}$ as a 
Farey fraction, where $r_{\nu} = h/k$, $0 < h \leq k$, $(h, k) = 1$, $k \leq N$, 
$1 \leq \nu \leq \varPhi(N)$, $\varPhi(N) = \varphi(1) + \dots + \varphi(N)$, $\varphi(\nu)$ 
the Euler function, Franel \cite{frj:fnp} and Landau \cite{lae:bhf} proved that the 
Riemann Hypothesis holds if and only if 
\begin{equation*}
\sum_{\nu = 1}^{\varPhi(N)} |\delta_{\nu}| = O(N^{{1 \over 2} + \epsilon})~,\end{equation*}
where $\delta_{\nu} = r_{\nu} - {\nu/\varPhi(N)}$, $\epsilon > 0$ arbitrary, is true. 
Intuitively, if $\nu / \varPhi(N)$ falls not far from the corresponding $r_{\nu}$, we could 
measure $\delta_{\nu}$ by $1/\varPhi(N)$. In \S \ref{min:ree}, we will prove 
\begin{rmmin}\label{fly:vld}\begin{equation}\label{fqu:vee}
\sum_{\nu = 1}^{\varPhi(N)} |\delta_{\nu}| = O(N^{1 \over 2} \log N)~.\end{equation}\end{rmmin}

Prior to Riemann's work, Dirichlet\footnote{Since his proof for the existence of 
infinite primes in a given arithmetic progression was (partly) based on his class 
number formula, Dirichlet assumed real variable $s$ for $L$-function. 
See~\cite[\S 1--6]{dah:mnt}.} discussed the distribution of primes in a given 
arithmetic progression $p \equiv a\!\pmod q$, $(a, q) = 1$, through $L$-function
\begin{equation*}L(s, \chi) = \sum_{n = 1}^{\infty} {\chi(n) \over n^{s}}~,\end{equation*}
where $\chi$ a Dirichlet character modulus $q$, $s$ a complex number, $n$ 
an integer, $Re(s) > 1$. Expanding $L$-function to the whole plane, for primitive 
character $\chi$ modulus $q$, we have functional equation\footnote{For its proof, 
see \cite[\S 9]{dah:mnt}.} analogous to \eqref{fuq:zte}
\begin{equation}\label{fuq:lte}
\!\!\!\! \biggl(\!{\pi \over q}\!\biggr)^{- {1 \over 2}(1 - s + \sigma_{\chi})} 
            \!\!\!\!\!\!\Gamma \biggl(\!\!\frac{1 - s + \sigma_{\chi}}{2}\!\!\biggr) L(s, \Bar{\chi}) = 
                                       \frac{i^{\sigma_{\chi}}q^{1 \over 2}}{\tau(\chi)}
                                                \biggl(\!\frac{\pi}{q}\!\biggr)^{- {1 \over 2}(s + \sigma_{\chi})}
                                                           \!\!\!\!\!\!\Gamma \biggl(\!\!\frac{s + \sigma_{\chi}}{2}\!\!\biggr) L(s, \chi),
\end{equation}
where $\tau(\chi) = \sum_{m = 1}^{q} \chi(m) \exp(2 \pi i m/q)$ the Gaussian sum, 
\[\sigma_{\chi} = \begin{cases}
                                 0  &\text{if $\chi(- 1) = 1$,} \\ 
                                 1  &\text{if $\chi(- 1) = - 1$.}\end{cases}\]
Through equation \eqref{fuq:lte}, we know that $L$-function has trivial zeros 
\[s = \begin{cases}
         -2, -4, -6, \ldots  &\text{when $\sigma_{\chi} = 0$,} \\
         -1, -3, -5, \ldots  &\text{when $\sigma_{\chi} = 1$.}\end{cases}\]Generalizing the 
Hypothesis, we would hope that all nontrivial zeros of $L(s, \chi)$ lie on line 
$Re(s) = 1/ 2$. In \S \ref{rst:gms}, we will prove 
\begin{lfmin}\label{lof:vld}All nontrivial zeros of $L(s, \chi)$ lie on line $Re(s) = 1/ 2$.\end{lfmin}
The validity of {\bf Theorem\/}~\ref{fly:vld} implies that\begin{gather}
\vartheta(x, \chi) = \sum_{n \leq x} \mu(n) \chi(n) = O(x^{1 \over 2} \log x)~,\label{tin:poo} 
\intertext{where $\chi$ a character modulus $q$, and} 
\frac{1}{L(s, \chi)} = \sum_{n = 1}^{\infty} \frac{\mu(n) \chi(n)}{n^{s}} 
                                  = s \int_{1}^{\infty} \frac{\vartheta(t, \chi)}{t^{s + 1}} dt~.\label{tin:wuo}
\end{gather}
where $Re(s) > 1$. Hence, we have {\bf Theorem\/}~\ref{lof:vld} by analytic continuation. 
\section{Theorem~\ref{fly:vld}}\label{min:ree}
\begin{hrlem}\label{hwn:thr}
\begin{equation*}\varPhi(N) = \varphi(1) + \dots + \varphi(N) = \frac{3{N^{2}}}{\pi^{2}} + O(N \log N)~. 
\end{equation*}\end{hrlem}
\begin{proof}This is Theorem 330 of Hardy and Wright \cite[p.~268]{hag:itn:wre}.\end{proof}
\begin{eulem}\label{eun:thr}
\begin{equation*} \sum_{n \leq x} \frac{1}{n} = \log x + \gamma + O\biggl(\frac{1}{x}\biggr)~,
\end{equation*}
where $\gamma$ is Euler constant.\end{eulem}
\begin{proof}This is Theorem 442 of Hardy and Wright \cite[p.~347]{hag:itn:wre}.\end{proof}
\begin{proof}[Proof of Theorem~\ref{fly:vld}]
For large $N$, denote $c_{j}$ constant, where $j \in \mathbb{Z}$,
\begin{gather*}\mathcal{F}_{N} = \{ h / k |~0 < h \leq k, k \leq N, (h, k) = 1 \}~, \\
\intertext{and} \mathcal{F}_{Nk} = \{ h/k |~h/k \in \mathcal{F}_{N}, k~\text{fixed} \}~.\end{gather*}

We build our proof based on following two {\em Observations\/}. In the first place, 
we look at instances when fraction $\nu / \varPhi(N)$ is {\em far away\/}, in terms of
$1 / \varPhi(N)$, from the corresponding $r_{\nu} \in \mathcal{F}_{Nk}$. By our instinct, 
if there is a restraint for such spot, we would have an easy job to evaluate the sum in 
\eqref{fqu:vee}. In the second place, we try to set factor $N^{1 \over 2}$ apart from the 
rest of terms so that we could get a desired outcome.
\begin{ovlap}\label{ovo:hun}For a fixed denominator $k$, we assume that there are 
$\lambda_{1}$ fractions $\nu_{\alpha} / \varPhi(N)$ corresponding to 
$r_{\nu_{\alpha}} \in \mathcal{F}_{Nk}$, where 
$\alpha = \alpha_{1}, \alpha_{2}, \dots, \alpha_{\lambda_{1}}$, 
fall into a $k$-interval $(j / k, (j + 1) / k)$, $1 \leq j \leq k - 1$;
there are $\lambda_{2}$ Farey fractions $r_{\nu_{\beta}} \in \mathcal{F}_{N}$ such that 
\begin{equation*}r_{\nu_{\beta}} \not \in \mathcal{F}_{Nk}, \text{and}~
r_{\nu_{\alpha_{1}}} \leq r_{\nu_{\beta}} \leq r_{\nu_{\alpha_{\lambda_{1}}}},\end{equation*}
where $\beta = \beta_{1}, \beta_{2}, \dots, \beta_{\lambda_{2}}$; 
there are $\lambda_{3}$ fractions $\nu_{\gamma} / \varPhi(N)$ corresponding to 
$r_{\nu_{\gamma}}$ such that $r_{\nu_{\gamma}} \not \in \mathcal{F}_{Nk}$, 
$\nu_{\alpha_{\lambda_{1}}} \leq \nu_{\gamma} \leq \nu_{\alpha_{\lambda_{1} + 1}}$, 
$\gamma = \gamma_{1}, \gamma_{2}, \ldots, \gamma_{\lambda_{3}}$. 
Since there is a one on one correspondence between $\nu / \varPhi(N)$ 
and Farey fraction $r_{\nu} \in \mathcal{F}_{Nk}$, we have
\begin{gather}\frac{\lambda_{1} + \lambda_{2} - 1}{\varPhi(N)} < \frac{1}{k}~, \label{sue:fst}\\ 
\frac{\lambda_{1} - 1}{k} < \frac{\lambda_{3} - 1}{\varPhi(N)}~, \label{sue:snd}\\
\intertext{and}(\lambda_{1} - 1) \cdot (\lambda_{1} + \lambda_{2} - 1) < \lambda_{3} - 1~.
\label{sue:trd}\end{gather}

Note that on average there are $\varPhi(N) / k - \varphi(k)/k$ Farey fractions 
$r_{\nu} \not \in \mathcal{F}_{Nk}$ fall in a $k$-interval $(j/k, (j + 1)/k)$, $1 \leq j \leq k$, 
the assumption tells us that $k$-intervals between $r_{\nu_{\alpha_{1}}}$ and 
$r_{\nu_{\alpha_{\lambda_{1}}}}$ are {\em void\/}, while $k$-interval 
$(r_{\nu_{\alpha_{\lambda_{1}}}}, r_{\nu_{\alpha_{\lambda_{1} + 1}}})$ is {\em abundant\/}
in terms of $r_{\nu} \not \in \mathcal{F}_{Nk}$. 
Since all $r_{\nu}^{\prime} \in \mathcal{F}_{N\iota}$, $\iota$ a prime, 
are equally distributed over $(0, 1)$, the {\em void\/} and the {\em abundant\/} are 
caused by $r_{\nu}^{\prime\prime} \in \mathcal{F}_{N\kappa}$, $\kappa$ a composite 
number. If $h/m \in \mathcal{F}_{N}$, $h/m \not \in \mathcal{F}_{Nk}$ falls into $k$-interval
$(r_{\nu_{\alpha_{\lambda_{1}}}}, r_{\nu_{\alpha_{\lambda_{1} + 1}}})$, but not previous
$\lambda_{1} - 1$ $k$-intervals, we have\begin{gather}
mr_{\nu_{\alpha_{\lambda_{1}}}} \leq h \leq mr_{\nu_{\alpha_{\lambda_{1} + 1}}}~, \notag \\ 
m (\lambda_{1} - 1) \leq k~. \label{ieq:kim}\end{gather}
\begin{siuon}\label{sin:koo}From equation~\eqref{ieq:kim}, we know that 
\begin{equation*}m < \sqrt{k}~, \text{if $\lambda_{1} > \sqrt{k} + 1$.}\end{equation*}
According to {\bf Lemma\/}~\ref{hwn:thr}, there are at most $O(k)$ such Farey fractions. 
Since $k \leq N$, we have\begin{equation*}
(\lambda_{1} - 1) \Biggl(\frac{\varPhi(N) - \varphi(k)}{k} 
                               + o\biggl(\frac{\varPhi(N) - \varphi(k)}{k} \biggr)\Biggr) 
                               - c_{1}k \geq \frac{\varPhi(N)}{k}~.\end{equation*}
Hence, $\lambda_{2} = O(\varPhi(N)/k)$, $\lambda_{1} = O(1)$ by~\eqref{sue:fst};\end{siuon}
\begin{siund}\label{sin:too}By equation~\eqref{ieq:kim}, we know that\begin{equation*}
m \leq k / (\lambda_{1} - 1)~, \text{if $\lambda_{1} \leq \sqrt{k} + 1$.}\end{equation*}
There are at most $O(k^{2} / (\lambda_{1} - 1)^{2})$ such Farey fractions. 
However, if $k \leq N^{2 \over 3}$, we have\begin{equation*}
(\lambda_{1} - 1) \Biggl(\frac{\varPhi(N) - \varphi(k)}{k} 
                               + o\biggl(\frac{\varPhi(N) - \varphi(k)}{k}\biggr)\Biggr)
                        - c_{2} \biggl(\frac{k^{2}}{(\lambda_{1} - 1)^{2}}\biggr)\!\! \geq \frac{\varPhi(N)}{k}~.
\end{equation*}Hence, $\lambda_{1} = O(1)$;\end{siund}
\begin{siurd}\label{sin:sai}If $\lambda_{1} \leq \sqrt{k} + 1$ and $k > N^{2 \over 3}$, we use
equation~\eqref{sue:trd}. Since\begin{align*}\lambda_{1} + \lambda_{2} - 1 &= 
(\lambda_{1} - 1)\Biggl(1 +(1 + o)\biggl(\frac{\varPhi(N) - \varphi(k)}{k}\biggr) 
                                - c_{3} \biggl(\frac{k^{2}}{(\lambda_{1} - 1)^{3}}\biggr)\Biggr), \\
\lambda_{3} - 1 &= (1 + o) \biggl(\frac{\varPhi(N) - \varphi(k)}{k}\biggr) + 
                                c_{4} \biggl(\frac{k^{2}}{(\lambda_{1} - 1)^{2}}\biggr),\end{align*}
we would get a contradiction if $\lambda_{1} = O(k^{\eta})$, where $0 < \eta \leq 1/2$. 
Therefore, $\lambda_{1} = O(1)$.\end{siurd}\end{ovlap}
\begin{totrm}For a fixed $k \leq N$, let each $k$-interval has equal number 
of Farey fractions, we see that fraction $\nu / \varPhi(N)$ would coincide with 
$r_{\nu} \in \mathcal{F}_{Nk}$. Hence, only uneven distribution of 
$r_{\nu} \not \in \mathcal{F}_{Nk}$ might cause $\nu / \varPhi(N)$ deviates from

$r_{\nu}$. We study following cases with consideration of 
{\em Observation\/}~\ref{ovo:hun}, and\begin{quotation}\begin{ounga}There should 
be no $k$-interval $(j/k, (j + 1)/k)$, $1 \leq j \leq k - 1$, $k \leq N$, with negative number 
of $r_{\nu} \not \in \mathcal{F}_{Nk}$, $r_{\nu} \in (j/k, (j + 1)/k)$.\end{ounga}\end{quotation}
\begin{kisei}If $k \leq N^{3 \over 4}$, then $k^{2} \leq N^{3 \over 2}$. On average, 
fraction $\nu / \varPhi(N)$ might deviate from the corresponding Farey 
fraction by $O(N^{3 \over 2} / k)$. By {\bf Lemma\/}~\ref{hwn:thr}, we know 
that\begin{equation*}
k \frac{1}{\varPhi(N)} \cdot c_{5} \cdot \frac{N^{3 \over 2}}{k} = \frac{c_{6}}{N^{1 \over 2}}
                                         = \frac{c_{6}N^{1 \over 2}}{N} \leq \frac{c_{7}N^{1 \over 2}}{k}~;
\end{equation*}\end{kisei}
\begin{kisbe}If $k > N^{3 \over 4}$, then on average fraction $\nu / \varPhi(N)$
might deviate from the corresponding Farey fraction by $O(k^{2} / k)$. Applying 
{\bf Lemma\/}~\ref{hwn:thr}, we know that\begin{equation}\label{pee:guu}
k \frac{1}{\varPhi(N)} \cdot c_{8} \cdot k = \frac{c_{9}k^{2}}{N^{2}}~.\end{equation}
There exist a number $\varpi_{1}(k)$ such that 
$1/\varpi_{1}(k) < k^{2}/N^{2} < N^{1 \over 2}/\varpi_{1}(k)$, i.e., 
$N^{2}/k^{2} < \varpi_{1}(k) < N^{5 \over 2}/k^{2}$. Substituting it into 
equation~\eqref{pee:guu}, we have\begin{equation*}
k \frac{1}{\varPhi(N)} \cdot c_{8} \cdot k \leq \frac{c_{10} N^{\frac{1}{2}}}{\varpi_{1}(k)}~;
\end{equation*}\end{kisbe}
\begin{kbead}In an extreme case, if fraction $\nu / \varPhi(N)$ deviates from $r_{\nu}$ 
by $O(k^{2})$, we have $k^{2} \leq O(\varPhi(N)/k)$ by {\em Condition N\/}. Therefore, 
$k \leq O(N^{2 \over 3})$, and\begin{equation}\label{pee:muu}
k \frac{1}{\varPhi(N)} \cdot c_{11} \cdot k^{2} = c_{12} \frac{k^{3}}{N^{2}}~.\end{equation}
There exist a number $\varpi_{2}(k)$ such that 
$1/\varpi_{2}(k) < k^{3}/N^{2} < N^{1 \over 2}/\varpi_{2}(k)$, i.e., 
$N^{2}/k^{3} < \varpi_{2}(k) < N^{5 \over 2}/k^{3}$. Substituting it into
equation~\eqref{pee:muu}, we have\begin{equation}\label{fna:muu}
k \frac{1}{\varPhi(N)} \cdot c_{11} \cdot k^{2} \leq \frac{c_{13} N^{\frac{1}{2}}}{\varpi_{2}(k)}~.
\end{equation}\end{kbead}\end{totrm}
In general, if there is a function $g(k)$ such that $\nu / \varPhi(N)$ deviates 
from the corresponding Farey fraction by $O(g(k)/k)$, through {\em Condition N\/} 
we have $g(k) \leq O(N^{2})$. Subsequently, we could find a corresponding 
$\varpi_{0}(k)$ such that $N^{2}/g(k) < \varpi_{0}(k) < N^{5 \over 2}/g(k)$, and get an 
inequality analogous to~\eqref{fna:muu}. Summing up these type of 
inequalities over $k$, we would have the desired consequence by 
{\bf Lemma\/}~\ref{eun:thr}.\end{proof}
\begin{comin}\label{cin:rhe}
All nontrivial zeros of $\zeta(s)$ lie on line $Re(s) = 1/ 2$.\end{comin}
\begin{proof}By Franel~\cite{frj:fnp}, Landau~\cite{lae:bhf} and 
{\bf Theorem~\ref{fly:vld}\/}. See also \cite[\S 12.2]{edh:rzf}.\end{proof}
\begin{comoo}\label{cin:moe}
\begin{equation*}M(x) = \sum_{n \leq x} \mu(n) = O(x^{1 \over 2}\log x)~.\end{equation*}
\end{comoo}
\begin{proof}By {\bf Theorem\/}~\ref{fly:vld}. See \cite[\S 12.2, p.~265]{edh:rzf}.\end{proof}
\section{Theorem~\ref{lof:vld}}\label{rst:gms}
\begin{palem}\label{sud:noo}Let $u(n)$ and $f(n)$ be arithmetic functions. Define 
sum function
\begin{equation*}U(t) = \sum_{n \leq t} u(n)~.\end{equation*}
Let $a$ and $b$ be nonnegtive integers with $a < b$. Then
\begin{equation*}
\!\!\!\sum_{n = a + 1}^{b}u(n)f(n) = U(b)f(b) - U(a)f(a +1) - \sum_{n = a + 1}^{b - 1}U(n)(f(n + 1) - f(n))~.
\end{equation*}
Let $x$ and $y$ be real numbers such that $0 \leq y < x$. If $f(t)$ is a function with a 
continuous derivative on the interval $[y, x]$, then
\begin{equation*}
\sum_{y < n \leq x} u(n)f(n) = U(x)f(x) - U(y)f(y) - \int_{y}^{x} U(t)f^{\prime}(t) dt~.\end{equation*}
\end{palem}
\begin{proof}This is Theorem A.~4 of Nathanson~\cite[p. 304]{nam:ant:tcb}.\end{proof}
\begin{tilem}\label{sud:roo}The estimation~\eqref{tin:poo} is valid.\end{tilem}
\begin{proof}For a principal character $\chi$, we use {\bf Corollary\/}~\ref{cin:moe}; for a 
nonprincipal character $\chi$, we use {\bf Lemma\/}~\ref{sud:noo}.\end{proof}
\begin{volem}\label{sud:woo}
The equation~\eqref{tin:wuo} is valid for $Re(s) > 1$.\end{volem}
\begin{proof}Through the Euler product\footnote{See \cite[\S 1]{dah:mnt}.} for $L(s, \chi)$
\begin{equation}\label{uqe:lrr}
L(s, \chi) = \prod_{p} \biggl(1 - \frac{\chi(p)}{p^{s}}\biggr)^{- 1}~.\end{equation}
where $p$ prime, $Re(s) > 1$, we have the first part of the equation~\eqref{tin:wuo}. 
For the second part, using {\bf Lemma\/}~\ref{sud:noo} we have
\begin{equation*}
\sum_{n = 1}^{x} \frac{\mu(n)\chi(n)}{n^{s}} = \frac{\vartheta(x, \chi)}{x^{s}} + 
                              s\int_{1}^{x} \frac{\vartheta(t, \chi)}{t^{s + 1}} dt~.\end{equation*}
Let $x \rightarrow \infty$, we get the conclusion.\end{proof}
\begin{proof}[Proof of Theorem~\ref{lof:vld}]
By {\bf Lemma\/}~\ref{sud:roo}, {\bf Lemma\/}~\ref{sud:woo} and analytic
continuation, we know that $1 / L(s, \chi)$ converges for all $Re(s) > 1/2$. Hence, 
$1 / L(s, \chi)$ is analytic for the half-plane $Re(s) > 1/2$. Through functional 
equation~\eqref{fuq:lte}, we see that $L(s, \chi)$ is symmetric with respect to 
$Re(s) = 1/2$ over $(0, 1)$. Therefore, all nontrivial zeros of $L(s, \chi)$ lie on 
line $Re(s) = 1/2$.\end{proof}
\end{document}